\documentclass[12pt,reqno,a4wide]{amsart}

\usepackage{tikz} 
\usepackage{calrsfs}
\usepackage{mathrsfs}
\usetikzlibrary{decorations.pathmorphing}
\usetikzlibrary{arrows.meta}
\usepgflibrary {shadings}
\usepackage{array}
\usepackage{xcolor}

\usepackage{comment}
\usepackage{hyperref}

\usetikzlibrary{decorations.pathreplacing}
\usetikzlibrary{fit}

\usepackage{pgfplots}

\allowdisplaybreaks

\catcode`,\active

\catcode`\,12

\begin{document}
\bibliographystyle{plain}

%
%

	\title[$k$-non-crossing trees and edge statistics modulo $k$]
	{$k$-non-crossing trees and edge statistics modulo $k$}

	\author[H. Prodinger ]{Helmut Prodinger }
	\address{Department of Mathematics, University of Stellenbosch 7602, Stellenbosch, South Africa
	and
NITheCS (National Institute for
Theoretical and Computational Sciences), South Africa.}
	\email{hproding@sun.ac.za}


	\begin{abstract}
		Instead of $k$-Dyck paths we consider the equivalent concept of $k$-non-crossing trees.
		This is our preferred approach relative to down-step statistics modulo $k$ (first studied by Heuberger, Selkirk, and Wagner by different methods). One symmetry argument
		about subtrees is needed and the rest goes along the lines of a paper by Flajolet and Noy.

	\end{abstract}
	

\maketitle

 
 \section{Non-crossing trees revisited}
 	
 	Assume that the nodes $1,\dots,n$ are arranged in a circle, call node 1 the root, and draw a tree using line segments such no crossings occur.
 	These objects are called \emph{non-crossing trees}. We only cite \cite{FN} and our own \cite{PP}, but there is much more literature that is not difficult to find.
 	Every node except for the root has two types of successors: left ones and right ones. See \cite{FN, PP}. Sometimes this is drawn as two trees that share a root node (`butterfly');
 	corresponding drawings are found in many papers on the subject. 
 	
\begin{figure} 	\label{f1}
\begin{tikzpicture}
	[every node/.style={circle,inner sep=2pt,minimum size=2em}]
	\def \n {10}
	\def \radius {3cm}
	\def \margin {6.5} 
	
	\foreach \s in {1,...,10}
	{
		\node(\s)[draw, circle,inner sep=2pt] at ({360/\n * (\s +1.48)}:\radius) {$\s$};
		arc ({360/\n * (\s - 1)+\margin}:{360/\n * (\s)-\margin}:\radius);
	}
	\draw[-]   (1) to (5);
	\draw[-]   (1) to (9);
	\draw[-]   (3) to (2);
	\draw[-]   (3) to (4);
	\draw[-]   (5) to (3);
	\draw[-]   (5) to (6);
	\draw[-]   (5) to (7);
	\draw[-]   (5) to (8);
	\draw[-]   (9) to (10);
\end{tikzpicture}
\begin{tikzpicture}[scale=0.6, every node/.style={circle,inner sep=3pt,minimum size=2em}]
	\node(1)[draw, circle,inner sep=3pt]  at (.0,.0){1};
	\node(5)[draw, circle,inner sep=3pt]  at (-2.0,-3.0){5};
	\node(9)[draw, circle,inner sep=3pt]  at (2.0,-3.0){9};
	\node(10)[draw, circle,inner sep=3pt  ]  at (4.0,-6.0){10};
	\node(3)[draw, circle,inner sep=3pt]  at (-4.0,-6.0){3};
	\node(4)[draw, circle,inner sep=3pt]  at (-2.0,-9.0){4};
	\node(2)[draw, circle,inner sep=3pt]  at (-6.0,-9.0){2};
	\node(6)[draw, circle,inner sep=3pt]  at (-1.0,-6.0){6};
	\node(7)[draw, circle,inner sep=3pt]  at (.5,-6.0){7};
	\node(8)[draw, circle,inner sep=3pt]  at (2.0,-6.0){8};

	\draw[-]   (1) to (5);
	\draw[-]   (1) to (9);
	\draw[-]   (3) to (2);
	\draw[-]   (3) to (4);
	\draw[-]   (5) to (3);
	\draw[-]   (5) to (6);
	\draw[-]   (5) to (7);
	\draw[-]   (5) to (8);
	\draw[-]   (9) to (10);
	\draw[->,>=stealth,ultra thick,red](-2.08,-3.8) to (5);
	\draw[->,>=stealth,ultra thick,red](2.015,-3.8) to (9);
	\draw[->,>=stealth,ultra thick,red](-4.035,-6.8) to (3);
\end{tikzpicture}
\caption{A non-crossing trees with 10 nodes and separators indicating where the non-root nodes split into the left part and the right part.}
\end{figure}
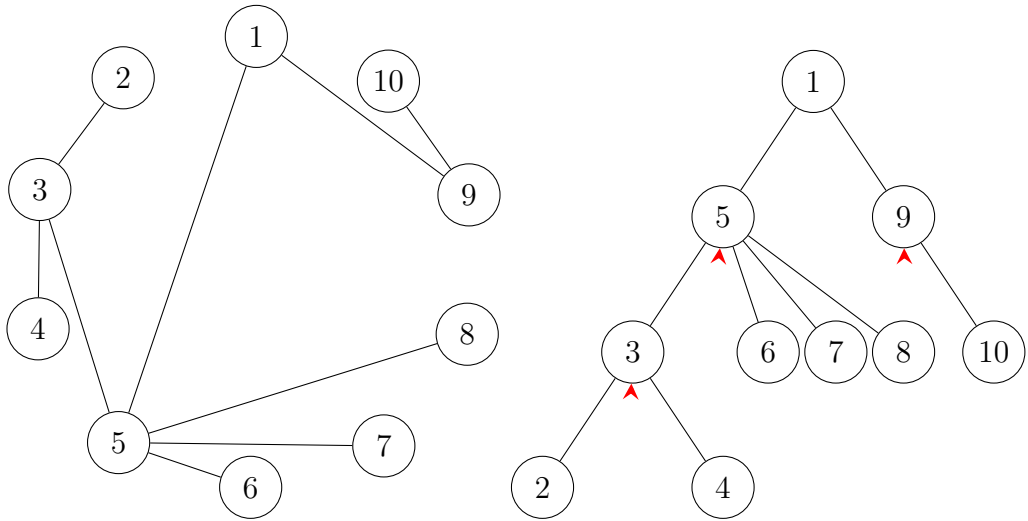

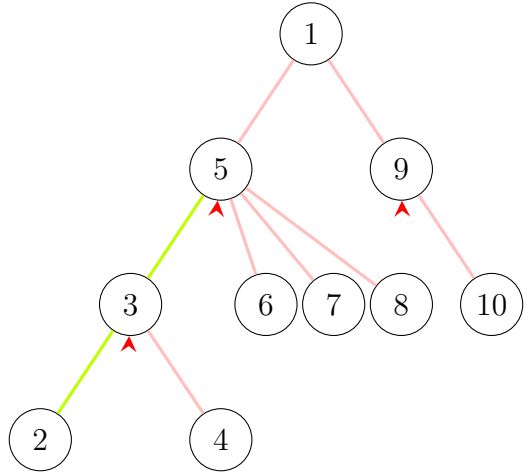
\begin{figure}
\begin{center}
\begin{tikzpicture}[scale=0.6, every node/.style={circle,inner sep=3pt,minimum size=2em}]
	\node(1)[draw, circle,inner sep=3pt]  at (.0,.0){1};
	\node(5)[draw, circle,inner sep=3pt]  at (-2.0,-3.0){5};
	\node(9)[draw, circle,inner sep=3pt]  at (2.0,-3.0){9};
	\node(10)[draw, circle,inner sep=3pt  ]  at (4.0,-6.0){10};
	\node(3)[draw, circle,inner sep=3pt]  at (-4.0,-6.0){3};
	\node(4)[draw, circle,inner sep=3pt]  at (-2.0,-9.0){4};
	\node(2)[draw, circle,inner sep=3pt]  at (-6.0,-9.0){2};
	\node(6)[draw, circle,inner sep=3pt]  at (-1.0,-6.0){6};
	\node(7)[draw, circle,inner sep=3pt]  at (.5,-6.0){7};
	\node(8)[draw, circle,inner sep=3pt]  at (2.0,-6.0){8};

	\draw[-,pink, very thick]   (1) to (5);
	\draw[-,pink, very thick]   (1) to (9);
	\draw[-,lime, very thick]   (3) to (2);
	\draw[-,pink, very thick]   (3) to (4);
	\draw[-,lime, very thick]   (5) to (3);
	\draw[-,pink, very thick]   (5) to (6);
	\draw[-,pink, very thick]   (5) to (7);
	\draw[-,pink, very thick]   (5) to (8);
	\draw[-,pink, very thick]   (9) to (10);
	\draw[->,>=stealth,ultra thick,red](-2.08,-3.8) to (5);
	\draw[->,>=stealth,ultra thick,red](2.015,-3.8) to (9);
	\draw[->,>=stealth,ultra thick,red](-4.035,-6.8) to (3);
\end{tikzpicture}

\end{center}
\caption{Left resp.\ right edges depicted in different colors; by design, the edges emanating from the root are all right edges.}
\end{figure}

It is interesting to note that \emph{even} trees as in \cite{C-even} are a similar concept to non-crossing trees.

We enumerate non-crossing trees with $n$ nodes and $j$ left and $n-1-j$ right edges. Clearly the total number of edges is $n-1$.
As one can see, the distribution isn't fair, as the root has all thes right edges as successors. We will use variables:
$n$ for the number of nodes and $\ell$ and $r$ for the two types of edges. We will use the butterfly decomposition due to Flajolet and Noy \cite{FN}.

\begin{equation*}
T=\frac z{1-B}, \quad B=\frac{T^2}{z};
\end{equation*}
$T$ stands for tree and $B$, which is only an auxiliary quantity, for butterfly. However, because of the anomaly of the root, we temporarily make $B$ the center of
interest:
\begin{equation*}
B=\frac{T^2}{z}=\frac{z}{(1-B)^2}.
\end{equation*}
Using the substitution $z=v(1-v)^2$, this can be solved, and the relevant solution is just $B=v$, and further $\frac Fz=\frac1{1-v}$. This can be extended with our extra variables
$\ell$ and $r$:
\begin{equation*}
	B=\frac{z}{(1-\ell B)(1-r B)},
\end{equation*}
then $z=v(1-\ell v)(1-rv)$ and $B=v$ and $\frac Fz=\frac{1}{1-rB}=\frac{1}{1-rv}$. Now we read off coefficients:
\begin{align*}
	[z^n]F&=[z^{n-1}]\frac Fz=[z^{n-1}]\frac1{1-rv}=\frac1{n-1}[z^{n-2}]\frac d{dz}\frac1{1-rv}\\
	&=\frac1{n-1}[z^{n-2}]\frac{dv}{dz}\frac d{dv}\frac1{1-rv}\\
	&=\frac1{n-1}[z^{n-2}]\frac1{1-2\ell v-2rv+3\ell rv^2} \frac{r}{(1-rv)^2}\\
	&=\frac1{n-1}\frac1{2\pi i}\oint\frac{dz}{z^{n-1}}\frac1{1-2\ell v-2rv+3\ell rv^2} \frac{r}{(1-rv)^2}\\
		&=\frac1{n-1}\frac1{2\pi i}\oint\frac{dv}{v^{n-1}(1-\ell v)^{n-1}(1-rv)^{n-1}}\frac{r}{(1-rv)^2}\\
				&=\frac r{n-1}[v^{n-2}]\frac{1}{(1-\ell v)^{n-1}(1-rv)^{n+1}}.
\end{align*}
As we can see, it is unnecessary to explicitly compute $\frac{dv}{dz}$ as it cancels out anyway. This will be very beneficial in the following sections. Furthermore,
\begin{align*}
	[z^n\ell^jr^{n-1-j}]F
	&=[v^{n-2}\ell^jr^{n-2-j}]\frac 1{n-1}\frac{1}{(1-\ell v)^{n-1}(1-rv)^{n+1}}\\
	&=\frac 1{n-1}\binom{n-2-j}{j}\binom{2n-2-j}{n-2-j}.
\end{align*}
This is the number of non-crossing trees with $n$ nodes, $j$ left edges and $n-1-j$ right edges.

\section{An application}

Lattice paths and certain types of trees are intimately related, and sometimes it is easier to analyze the trees instead of the
paths, an example being \cite{CMcL, prodinger-hills}. This will also happen here, as we will use the analysis of non-crossing trees from the Introduction to lattice paths.

We transform non-crossing trees into so-called $2$-Dyck paths: Up-steps $(1,1)$ are as usual, but there are down-steps
$(1,-2)$ of two units. Otherwise, the path must be non-negative and eventually return to the $x$-axis. 
For this transformation, we walk around the tree and translate down-steps into up-steps and vice versa. However, we need extra up-steps to keep  
the balance. For that we use the separators, and also draw them for end-nodes, so that there are $n-1$ such separator markers
present. Then, whenever we meet one, we also make an up-step. 
	
\begin{figure}
\begin{center}
\begin{tikzpicture}[scale=0.6, every node/.style={circle,inner sep=3pt,minimum size=2em}]
	\node(1)[draw, circle,inner sep=3pt]  at (.0,.0){1};
	\node(5)[draw, circle,inner sep=3pt]  at (-2.0,-3.0){5};
	\node(9)[draw, circle,inner sep=3pt]  at (2.0,-3.0){9};
	\node(10)[draw, circle,inner sep=3pt  ]  at (4.0,-6.0){10};
	\node(3)[draw, circle,inner sep=3pt]  at (-4.0,-6.0){3};
	\node(4)[draw, circle,inner sep=3pt]  at (-2.0,-9.0){4};
	\node(2)[draw, circle,inner sep=3pt]  at (-6.0,-9.0){2};
	\node(6)[draw, circle,inner sep=3pt]  at (-1.0,-6.0){6};
	\node(7)[draw, circle,inner sep=3pt]  at (-5.5,-12.0){7};
	\node(8)[draw, circle,inner sep=3pt]  at (-3.5,-12.0){8};

	\draw[-,pink, very thick]   (1) to (5);
	\draw[-,pink, very thick]   (1) to (9);
	\draw[-,lime, very thick]   (3) to (2);
	\draw[-,pink, very thick]   (3) to (4);
	\draw[-,lime, very thick]   (5) to (3);
	\draw[-,lime, very thick]   (4) to (7);
	\draw[-,pink, very thick]   (5) to (6);
	\draw[-,lime, very thick]   (4) to (8);
	\draw[-,pink, very thick]   (9) to (10);
	\draw[->,>=stealth,ultra thick,red](-2.08,-3.8) to (5);
	\draw[->,>=stealth,ultra thick,red](2.015,-3.8) to (9);
	\draw[->,>=stealth,ultra thick,red](-4.035,-6.8) to (3);
	\draw[->,>=stealth,ultra thick,red](-6.00,-9.8) to (2);
	\draw[->,>=stealth,ultra thick,red](-2.00,-9.8) to (4);
	\draw[->,>=stealth,ultra thick,red](-1.00,-6.8) to (6);
	\draw[->,>=stealth,ultra thick,red](-5.5,-12.8) to (7);
	\draw[->,>=stealth,ultra thick,red](-3.5,-12.8) to (8);
	\draw[->,>=stealth,ultra thick,red](4.00,-6.8) to (10);
	
\end{tikzpicture}
\end{center}\vskip0.5cm
\begin{center}
\begin{tikzpicture}[scale=0.4 ]\label{brown}

	\foreach \i in {0,...,8} {
		\draw [very thin,gray] (0,\i) to (27,\i) ;   
	}
	\foreach \i in {0,...,27} {
		\draw [very thin,gray] (\i,0) to (\i,8) ;   
	}
	
	 \draw[thick](0,0) to (4,4);
	 	 \draw[thick, blue](4,4) to (5,2);
	 	 	 \draw[thick](5,2) to (9,6);
	 	 	 	 	 	 \draw[very thick,blue](9,6) to (10,4);
	 	 	 	 	 	 	 	 	 \draw[thick](10,4) to (12,6);
	 	 	 	 	 	 	 	 	 \draw[very thick,blue](12,6) to (13,4);
	 	 	 	 	 	 	 	 	 \draw[thick](13,4) to (14,5);
	 	 	 	 	 	 	 	 	 	 	 	 \draw[very thick,brown](14,5) to (16,1);
	 	 	 	 	 	 	 	 	 	 	 	 	 	 	 \draw[thick](16,1) to (19,4);
	 	 	 	 	 	 	 	 	 	 	 	 	 	 	 	 	 	 	 	\draw[very thick,blue](19,4) to (21,0);
	 	 	 	 	 	 	 	 	 	 	 	 	 	 	 	 	 	 	 	 	 	 	 	 	 	 	 \draw[thick](21,0) to (25,4);
	 	 	 \draw[very thick,blue](25,4) to (27,0);
\end{tikzpicture}
\end{center}
\caption{A non-crossing tree and the corresponding 2-Dyck path.}
\end{figure}
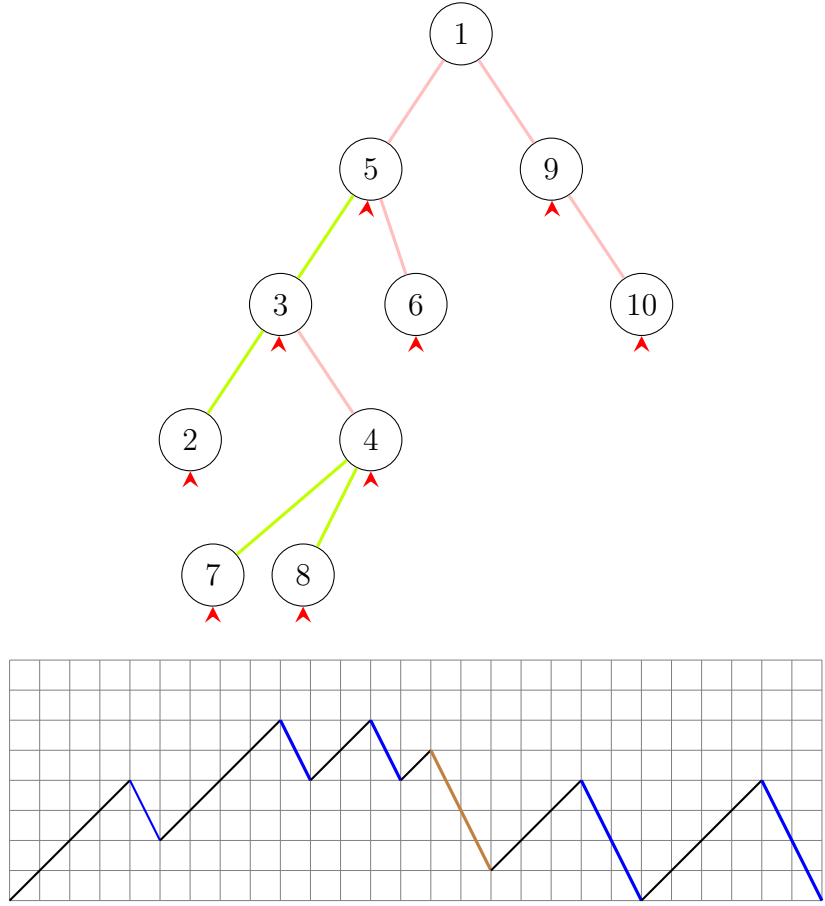

In the example we get the path in Figure~{\ref{brown}}.

 Note that $n$ nodes of the tree correspond to $n-1$ down-steps. More of such considerations can be found in
     \cite{Selkirk-master, GPW, prodinger-jcmcc}.

\begin{figure}

\begin{center}
	\begin{tikzpicture}[scale=0.6, every node/.style={circle,inner sep=3pt,minimum size=2em}]
		\node(1)[draw, circle,inner sep=3pt]  at (.0,.0){1};
		\node(5)[draw, circle,inner sep=3pt]  at (-2.0,-3.0){5};
		\node(9)[draw, circle,inner sep=3pt]  at (2.0,-3.0){9};
		\node(10)[draw, circle,inner sep=3pt  ]  at (4.0,-6.0){10};
		\node(3)[draw, circle,inner sep=3pt]  at (-4.0,-6.0){3};
		\node(2)[draw, circle,inner sep=3pt]  at (-2.0,-9.0){2};
		\node(4)[draw, circle,inner sep=3pt]  at (-6.0,-9.0){4};
		\node(6)[draw, circle,inner sep=3pt]  at (-1.0,-6.0){6};
		\node(7)[draw, circle,inner sep=3pt]  at (-9.5,-12.0){7};
		\node(8)[draw, circle,inner sep=3pt]  at (-7.5,-12.0){8};

		\draw[-,pink, very thick]   (1) to (5);
		\draw[-,pink, very thick]   (1) to (9);
		\draw[-,pink, very thick]   (3) to (2);
		\draw[-,lime, very thick]   (3) to (4);
		\draw[-,lime, very thick]   (5) to (3);
		\draw[-,lime, very thick]   (4) to (7);
		\draw[-,pink, very thick]   (5) to (6);
		\draw[-,lime, very thick]   (4) to (8);
		\draw[-,pink, very thick]   (9) to (10);
		\draw[->,>=stealth,ultra thick,red](-2.08,-3.8) to (5);
		\draw[->,>=stealth,ultra thick,red](2.015,-3.8) to (9);
		\draw[->,>=stealth,ultra thick,red](-4.035,-6.8) to (3);
		\draw[->,>=stealth,ultra thick,red](-2.00,-9.8) to (2);
		\draw[->,>=stealth,ultra thick,red](-6.00,-9.8) to (4);
		\draw[->,>=stealth,ultra thick,red](-1.00,-6.8) to (6);
		\draw[->,>=stealth,ultra thick,red](-9.5,-12.8) to (7);
		\draw[->,>=stealth,ultra thick,red](-7.5,-12.8) to (8);
		\draw[->,>=stealth,ultra thick,red](4.00,-6.8) to (10);
	\end{tikzpicture}%
	\begin{tikzpicture}[scale=0.6, every node/.style={circle,inner sep=3pt,minimum size=2em}]
		\node(1)[draw, circle,inner sep=3pt]  at (.0,.0){1};
		\node(5)[draw, circle,inner sep=3pt]  at (-2.0,-3.0){5};
		\node(9)[draw, circle,inner sep=3pt]  at (2.0,-3.0){9};
		\node(10)[draw, circle,inner sep=3pt  ]  at (4.0,-6.0){10};
		\node(3)[draw, circle,inner sep=3pt]  at (-4.0,-6.0){3};
		\node(2)[draw, circle,inner sep=3pt]  at (-2.0,-9.0){2};
		\node(4)[draw, circle,inner sep=3pt]  at (-6.0,-9.0){4};
		\node(6)[draw, circle,inner sep=3pt]  at (-1.0,-6.0){6};
		\node(7)[draw, circle,inner sep=3pt]  at (-4.5,-12.0){7};
		\node(8)[draw, circle,inner sep=3pt]  at (-1.5,-12.0){8};

		\draw[-,pink, very thick]   (1) to (5);
		\draw[-,pink, very thick]   (1) to (9);
		\draw[-,pink, very thick]   (3) to (2);
		\draw[-,lime, very thick]   (3) to (4);
		\draw[-,lime, very thick]   (5) to (3);
		\draw[-,pink, very thick]   (4) to (7);
		\draw[-,pink, very thick]   (5) to (6);
		\draw[-,pink, very thick]   (4) to (8);
		\draw[-,pink, very thick]   (9) to (10);
		\draw[->,>=stealth,ultra thick,red](-2.08,-3.8) to (5);
		\draw[->,>=stealth,ultra thick,red](2.015,-3.8) to (9);
		\draw[->,>=stealth,ultra thick,red](-4.035,-6.8) to (3);
		\draw[->,>=stealth,ultra thick,red](-2.00,-9.8) to (2);
		\draw[->,>=stealth,ultra thick,red](-6.00,-9.8) to (4);
		\draw[->,>=stealth,ultra thick,red](-1.00,-6.8) to (6);
		\draw[->,>=stealth,ultra thick,red](-4.5,-12.8) to (7);
		\draw[->,>=stealth,ultra thick,red](-1.5,-12.8) to (8);
		\draw[->,>=stealth,ultra thick,red](4.00,-6.8) to (10);
	\end{tikzpicture}
\end{center}
\caption{Transforming the tree. The number of green edges corresponds to the brown down-steps.}
\end{figure}
The goal is to match the brown down-steps to the left edges, say.
In particular, the interest is, on which level modulo $k$ they land (or, equivalently, start). First, the tree needs to be modified. The reason is this decomposition in Figure~\ref{deco}. Indeed, $T_1$ ``sits'' on level 1 (odd) but a subtree of $T_1$ ``sits'' on level 2 (even). So we need to swap subtrees in such a case. The next section will provide more details.
\begin{figure}[h]
	\begin{center}
		\begin{tikzpicture}[scale=0.4]
			
			(20,6); \draw[ultra thick] (0.0,0.) to (1.,1.);
			
			\draw (1,1) .. controls (2,5) .. (3,1); \draw
			(3,1) .. controls (4,7) .. (5,1); \node at
			(5.5,1) {$\cdot$}; \node at (6,1) {$\cdot$};
			\node at (6.5,1) {$\cdot$}; \draw (7,1) ..
			controls (8,3) .. (9,1); \draw [ultra
			thick](9,1) to  (10,2); \draw (1+9,1+1) ..
			controls (2+9,6) .. (3+9,1+1); \draw (3+9,1+1)
			.. controls (4+9,4) .. (5+9,1+1); \node at
			(5.5+9,2) {$\cdot$}; \node at (6+9,2)
			{$\cdot$}; \node at (6.5+9,2) {$\cdot$}; \draw
			(16,2) .. controls (17,5) .. (18,2); \draw
			[ultra thick](18,2) to  (19,0);
			
			\node at (6.5+14,0) {$\cdots$};

			\draw [thick, decorate, decoration={brace,
				amplitude=10pt, mirror, raise=4pt}] (1cm, -0.5)
			to node[below,yshift=-0.5cm] {$T_{1}$} (9cm,
			-0.5);
			
			\draw [thick, decorate, decoration={brace,
				amplitude=10pt, mirror, raise=4pt}] (10cm,
			-0.5) to node[below,yshift=-0.5cm] {$T_{2}$}
			(18cm,- 0.5);

		\end{tikzpicture} \end    {center}
		
		\caption{The decomposition of  a 2-Dyck
			path.}
		\label{deco}
		\end{figure}

\section{Generalization}

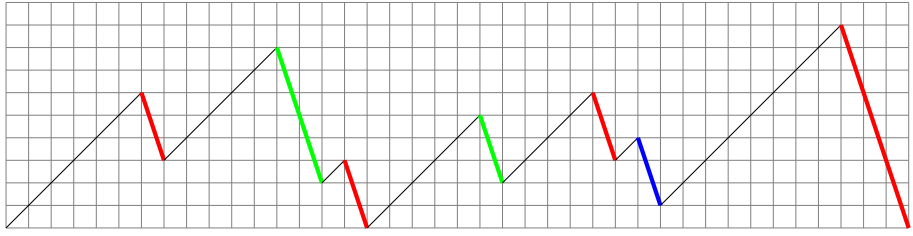
\begin{figure}
	 \label{threeDyck}

\begin{center}
	\begin{tikzpicture}[scale=0.3 ]

		\foreach \i in {0,...,10} {
			\draw [very thin,gray] (0,\i) to (40,\i) ;   
		}
		\foreach \i in {0,...,40} {
			\draw [very thin,gray] (\i,0) to (\i,10) ;   
		}
		
		\draw(0,0) to (6,6);
		\draw [ultra thick, red](6,6)-- (7,3);
		\draw(7,3) to (12,8);
		\draw [ultra thick, green](12,8)-- (14,2);
		\draw(14,2) to (15,3);
		\draw [ultra thick, red](15,3)-- (16,0);
		\draw(16,0) to (21,5);
		\draw [ultra thick, green](21,5)-- (22,2);
		\draw(22,2) to (26,6);
		\draw [ultra thick, red](26,6)-- (27,3);
		\draw(27,3) to (28,4);
		\draw [ultra thick, blue](28,4)-- (29,1);
		\draw(29,1) to (37,9);
		\draw [ultra thick, red](37,9)-- (40,0);
	\end{tikzpicture}
\end{center}
\caption{A 3-Dyck path. 6 down-steps land on level 0~$(\bmod\;3)$,  1 on level 1~$(\bmod\;3)$, and 3 on level 2~$(\bmod\;3)$}
\end{figure}

Instead of down-steps of two units and one separator, this works as well for down-steps $(1,-k)$ and $k-1$ separators.
Here (Figure~\ref{threeDyck}) is a 3-Dyck paths: 6 down-steps land on level 0~$(\bmod\;3)$,  1 on level 1~$(\bmod\;3)$, and 3 on level 2~$(\bmod\;3)$.

The butterfly equation is 
\begin{equation*}
B=\frac{T^k}{z}=\frac z{(1-r_1B)\dots (1-r_kB)}
\end{equation*}
with variables $r_1,\dots, r_k$ to count the down-steps ending (or beginning) on a level $\equiv i (\bmod k)$. Eventually $T=\frac{z}{1- r_kB}$. For the solution,
the substitution $z=v(1-r_1v)\dots(1-r_kv)$ works, and $B=v$, and thus $\frac Tz=\frac{1}{1-r_kv}$.
	
	Reading off coefficients is similar to the previous case $k=2$:
	\begin{align*}
[z^n]F&=[z^{n-1}]\frac Fz=\frac1{n-1}[z^{n-2}]\frac{d}{dz}\frac Fz\\
&=\frac1{n-1}[z^{n-2}]\frac{dv}{dz}\frac{d}{dv}\frac Fz=\frac1{n-1}[z^{n-2}]\frac{dv}{dz}\frac{r_k}{(1-r_kv)^2}\\
&=\frac1{n-1}\frac1{2\pi i}\oint \frac{dz}{z^{n-1}}\frac{dv}{dz}\frac{r_k}{(1-r_kv)^2}\\
&=\frac1{n-1}\frac1{2\pi i}\oint \frac{dv}{v^{n-1}(1-r_1v)^{n-1}\dots(1-r_kv)^{n-1}}\frac{r_k}{(1-r_kv)^2}\\
&=\frac1{n-1}[v^{n-2}]  \frac{r_k}{(1-r_1v)^{n-1}\dots(1-r_{k-1}v)^{n-1}(1-r_{k}v)^{n+1}}.
	\end{align*}
Furthermore (with $a_1+\dots +a_k=n-1$)
\begin{align*}
[z^nr_1^{a_1}\dots r_k^{a_k}]F&=\frac1{n-1}[v^{n-2}r_1^{a_1}\dots r_k^{a_k}]  \frac{r_k}{(1-r_1v)^{n-1}\dots(1-r_{k-1}v)^{n-1}(1-r_{k}v)^{n+1}}\\
&=\frac1{n-1}[v^{n-2}r_1^{a_1}\dots r_k^{a_k-1}]  \frac{1}{(1-r_1v)^{n-1}\dots(1-r_{k-1}v)^{n-1}(1-r_{k}v)^{n+1}}\\
&=\frac1{n-1}\binom{n-2+a_1}{a_1}\dots \binom{n-2+a_{k-1}}{a_{k-1}}\binom{n-1+a_k}{a_k-1}.
\end{align*}

This is the formula in Theorem~6 in \cite{HSW}, for $t=0$, and $n\to n+1$. For more general $-k<-t\le 0\ (\Leftrightarrow 0\le t<k)$, we will work this out in the next section.

The modification of the tree (rotation of the subtrees, from top to bottom, so that the down-step enumeration matches the edge enumeration) is as follows (for $k=3$)

	
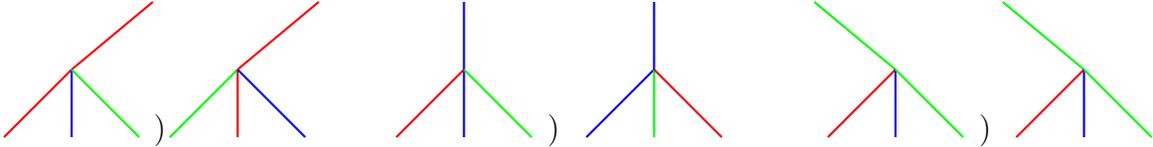
\begin{tikzpicture}
		[scale=0.6,thick,level 1/.style={sibling distance=15mm},
		level 2/.style={sibling distance=15mm}]
		\coordinate
		child[red,xshift=-1.8cm]   {child[red] child[blue] child[green]}
		;
	\end{tikzpicture}%
	\begin{tikzpicture} 
 \raisebox{0.8cm}{$\Rightarrow$};
\end{tikzpicture}\hspace{0.2cm}
\begin{tikzpicture} 
	[scale=0.6,thick,level 1/.style={sibling distance=15mm},
	level 2/.style={sibling distance=15mm}]
	\coordinate
	child[red,xshift=-1.8cm]   { child[green] child[red] child[blue]}
	;
\end{tikzpicture}\hspace{1cm}%
\begin{tikzpicture}
	[scale=0.6,thick,level 1/.style={sibling distance=15mm},
	level 2/.style={sibling distance=15mm}]
	\coordinate
	child[blue,xshift=0cm]   {child[red] child[blue] child[green]}
	;
\end{tikzpicture}\hspace{0.2cm}
\begin{tikzpicture} 
	\raisebox{0.8cm}{$\Rightarrow$};
\end{tikzpicture}\hspace{0.5cm}%
\begin{tikzpicture}
	[scale=0.6,thick,level 1/.style={sibling distance=15mm},
	level 2/.style={sibling distance=15mm}]
	\coordinate
	child[blue,xshift=0cm]   {child[blue] child[green] child[red] }
	;
\end{tikzpicture}\hspace{1.2cm}%
\begin{tikzpicture}
	[scale=0.6,thick,level 1/.style={sibling distance=15mm},
	level 2/.style={sibling distance=15mm}]
	\coordinate
	child[green,xshift=1.8cm]  {child[red] child[blue] child[green]}
	;
\end{tikzpicture}\hspace{0.2cm}
\begin{tikzpicture} 
	\raisebox{0.8cm}{$\Rightarrow$};
\end{tikzpicture}\hspace{0.3cm}%
\begin{tikzpicture}
	[scale=0.6,thick,level 1/.style={sibling distance=15mm},
	level 2/.style={sibling distance=15mm}]
	\coordinate
	child[green,xshift=1.8cm]   {child[red] child[blue] child[green]}
	;
\end{tikzpicture}

It is easy to figure out how this works more generally for $k$ successors instead of 3. It is always a cyclic shift, by $k-1,k-2,\dots,0$ positions, depending on the edge 
we are considering. 

\begin{figure}
  	
 \begin{center}
 	\begin{tikzpicture}[scale=0.35 ]

 		\foreach \i in {0,...,5} {
 			\draw [very thin,gray,dotted] (-1,\i) to (21,\i) ;   
 		}
 		\foreach \i in {-1,...,21} {
 			\draw [very thin,gray,dotted] (\i,0) to (\i,5) ;   
 		}
 		\draw [  thick ] (-1,1) to (21,1) ;  
 		\draw[  thick ](0,1) to (3,4);
 		 \draw[  thick ]  (3,4)--(4,2);
 		 \draw[  thick ]  (4,2)--(5,3);
 		 \draw[  thick ]  (5,3)--(6,1);
 		 \draw[  thick ]  (6,1)--(7,2);
\draw[  thick ]  (7,2)--(8,0);
\draw[  thick ]  (8,0)--(10,2);
\draw[  thick ]  (10,2)--(11,0);
\draw[  thick ,red]  (11,0)--(12,1);
\draw[  thick ]  (12,1)--(15,4);
\draw[  thick ]  (15,4)--(16,2);
\draw[  thick ]  (16,2)--(17,3);
\draw[  thick ]  (17,3)--(18,1);
\draw[  thick ]  (18,1)--(20,3);
\draw[  thick ]  (20,3)--(21,1);
 	\end{tikzpicture}\hspace*{0.5cm}%
 	 		\begin{tikzpicture}[scale=0.35 ]

 			\foreach \i in {0,...,5} {
 				\draw [very thin,gray,dotted] (-1,\i) to (21,\i) ;   
 			}
 			\foreach \i in {-1,...,21} {
 				\draw [very thin,gray,dotted] (\i,0) to (\i,5) ;   
 			}
 			\draw [  thick ] (-1,1) to (21,1) ;  
 			\draw[  thick ](0,1) to (3,4);
 			\draw[  thick ]  (3,4)--(4,2);
 			\draw[  thick ]  (4,2)--(5,3);
 			\draw[  thick ]  (5,3)--(6,1);
 			\draw[  thick ]  (6,1)--(7,2);
 			\draw[  thick ]  (7,2)--(8,0);
 			\draw[  thick ]  (8,0)--(10,2);
 			\draw[  thick ]  (10,2)--(11,0);
 			\draw[  thick ,red]  (-1,0)--(0,1);
 			\draw[  thick ]  (12,1)--(15,4);
 			\draw[  thick ]  (15,4)--(16,2);
 			\draw[  thick ]  (16,2)--(17,3);
 			\draw[  thick ]  (17,3)--(18,1);
 			\draw[  thick ]  (18,1)--(20,3);
 			\draw[  thick ]  (20,3)--(21,1);
 			 
 		\end{tikzpicture}
 		 \end{center}
 		   \caption{Decomposition of paths bounded by the line $y=-1$ into   two paths.}
 		   \label{blah}
 		   \end{figure}
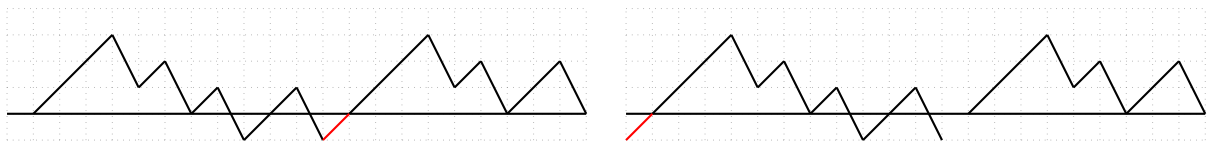

The following differentiation will be used in the sequel. It is just the differentiation of a product, as usual.

\begin{equation*}
\frac d{dv}\prod_{j=k-t}^k\frac1{1-r_jv}=\prod_{j=k-t}^k\frac1{1-r_jv}\cdot\sum_{i=k-t}^k\frac{r_i}{1-r_iv}.
\end{equation*}
Our strategy is to bijectively map $k$-Dyck paths bounded below by $y=-t$ into $t+1$ $k$-non-crossing trees of altogether
$n-t-1$ edges, and the special symbol attached to the root varies from $r_{k}$, $r_{k-1}$ \dots to $r_{k-t}$. Figure~\ref{blah} shows a small example, and more examples are in \cite{prodinger-jcmcc}.
Then
\begin{align*}
	[z^n]F&=[z^{n-t-1}]\frac F{z^{t+1}}=[z^{n-t-1}]\prod_{j=k-t}^k\frac1{1-r_jv}
	=\frac1{n-t-1}[z^{n-t-2}]\frac{d}{dz}\prod_{j=k-t}^k\frac1{1-r_jv}\\
	&=\frac1{n-t-1}[z^{n-t-2}]\frac{dv}{dz}\prod_{j=k-t}^k\frac1{1-r_jv}\cdot\sum_{i=k-t}^k\frac{r_i}{1-r_iv}\\
	&=\frac1{n-t-1}\frac1{2\pi i}\oint\frac{dz}{z^{n-t-1}}
	\frac{dv}{dz}\prod_{j=k-t}^k\frac1{1-r_jv}\cdot\sum_{i=k-t}^k\frac{r_i}{1-r_iv}\\
	&=\frac1{n-t-1}\frac1{2\pi i}\oint\frac{dv}{v^{n-t-1}(1-r_1v)^{n-t-1}\dots(1-r_kv)^{n-t-1}}
	 \prod_{j=k-t}^k\frac1{1-r_jv}\cdot\sum_{i=k-t}^k\frac{r_i}{1-r_iv}\\
	 &=\frac1{n-t-1}\frac1{2\pi i}\oint\frac{dv}{v^{n-t-1}\prod_{h=1}^{k-t-1}(1-r_hv)^{n-t-1} 
	 	\prod_{\ell=k-t}^{k}(1-r_\ell v)^{n-t}	 }
	 \sum_{i=k-t}^k\frac{r_i}{1-r_iv}\\
	 	 &=\frac1{n-t-1}[v^{n-t-2}] \frac{1}{ \prod_{h=1}^{k-t-1}(1-r_hv)^{n-t-1} 
	 	\prod_{\ell=k-t}^{k}(1-r_\ell v)^{n-t}	 }
	 \sum_{i=k-t}^k\frac{r_i}{1-r_iv}.
	 	 	\end{align*}
	 	 	%
	Furthermore ($a_1+\dots +a_k=n-t-1$)
	\begin{align*}
		[z^nr_1^{a_1}&\dots r_k^{a_k}]F\\& = \frac1{n-t-1}[v^{n-t-2}r_1^{a_1}\dots r_k^{a_k}] \sum_{i=k-t}^k\frac{r_i}{ \prod_{h=1}^{k-t-1}(1-r_hv)^{n-t-1} 
			\prod_{\ell=k-t}^{k}(1-r_\ell v)^{n+[i=\ell]-t}	 }\\
			&= \frac1{n-t-1}\sum_{i=k-t}^k \prod_{h=1}^{k-t-1} \binom{n-t-2+a_h}{a_h}\prod_{\ell=k-t}^{k}\binom{n+[i=\ell]-t-1+a_\ell-[i=\ell]}{a_\ell-[i=\ell]}  \\
			&= \frac1{n-t-1}\sum_{i=k-t}^k \prod_{h=1}^{k-t-1} \binom{n-t-2+a_h}{a_h}\prod_{\ell=k-t}^{k}\binom{n -t-1+a_\ell }{a_\ell-[i=\ell]}  \\
			&= \frac1{n-t-1}\sum_{i=k-t}^k \frac{a_i}{n-t}\prod_{h=1}^{k-t-1} \binom{n-t-2+a_h}{a_h}\prod_{\ell=k-t}^{k}\binom{n -t-1+a_\ell }{a_\ell }.
								\end{align*}
								The formula looks better when $n-t-1=N$; then it compares and matches with the formula from
								\cite{HSW}
					\begin{align*}
						 						 \frac{a_{k-t}+\dots+a_k}{N(N+1)}\prod_{h=1}^{k-t-1} \binom{N-1+a_h}{a_h}\prod_{\ell=k-t}^{k}\binom{ N+a_\ell }{a_\ell}. 
					\end{align*}
  
 \bibliographystyle{plain}


\begin{thebibliography}{1}
  	
  	\bibitem{C-even}
  	Naiomi~T. Cameron.
  	\newblock The combinatorics of even trees.
  	\newblock In {\em Proceedings of the {T}hirty-first {S}outheastern
  		{I}nternational {C}onference on {C}ombinatorics, {G}raph {T}heory and
  		{C}omputing ({B}oca {R}aton, {FL}, 2000)}, volume 147, pages 129--143, 2000.
  	
  	\bibitem{CMcL}
  	Naiomi~T. Cameron and Jillian~E. McLeod.
  	\newblock Returns and hills on generalized {D}yck paths.
  	\newblock {\em J. Integer Seq.}, 19(6):Article 16.6.1, 28, 2016.
  	
  	\bibitem{FN}
  	Philippe Flajolet and Marc Noy.
  	\newblock Analytic combinatorics of non-crossing configurations.
  	\newblock {\em Discrete Math.}, 204(1-3):203--229, 1999.
  	
  	\bibitem{GPW}
  	Nancy S.~S. Gu, Helmut Prodinger, and Stephan Wagner.
  	\newblock Bijections for a class of labeled plane trees.
  	\newblock {\em European J. Combin.}, 31(3):720--732, 2010.
  	
  	\bibitem{HSW}
  	Clemens Heuberger, Sarah~J. Selkirk, and Stephan Wagner.
  	\newblock Enumeration of generalized {D}yck paths based on the height of
  	down-steps modulo {$k$}.
  	\newblock {\em Electron. J. Combin.}, 30(1):Paper No. 1.26, 18, 2023.
  	
  	\bibitem{PP}
  	Alois Panholzer and Helmut Prodinger.
  	\newblock Bijections for ternary trees and non-crossing trees.
  	\newblock {\em Discrete Math.}, 250(1-3):181--195, 2002.
  	
  	\bibitem{prodinger-jcmcc}
  	Helmut Prodinger.
  	\newblock On $k$-{D}yck paths with a negative boundary.
  	\newblock {\em Journal of Combinatorial Mathematics and Combinatorial
  		Computing}.
  	
  	\bibitem{prodinger-hills}
  	Helmut Prodinger.
  	\newblock Returns, hills, and {$t$}-ary trees.
  	\newblock {\em J. Integer Seq.}, 19(7):Article 16.7.2, 8, 2016.
  	
  	\bibitem{Selkirk-master}
  	Sarah~J. Selkirk.
  	\newblock Msc-thesis: On a generalisation of $k$-{D}yck paths.
  	\newblock {\em Stellenbosch University}, 2019.
  	
  \end{thebibliography}

\end{document}